%% file: m7-5.tex
\documentclass{gtart_h}

\input gtmonout

\volumenumber{7}
\volumename{Proceedings of the Casson Fest}
\volumeyear{2004}
\papernumber{5}
\pagenumbers{135}{144}
\received{19 March 2004}
\revised{28 July 2004}
\accepted{4 August 2004}
\published{18 September 2004}

\usepackage{amsmath,amssymb,rlepsf}
\nocolon

\begin{document}
\newtheorem{theorem}{Theorem}
\newtheorem{prop}[theorem]{Proposition}
\newtheorem{lemma}[theorem]{Lemma}
\newtheorem{claim}[theorem]{Claim}
\newtheorem{conj}[theorem]{Conjecture}
\newtheorem{cor}[theorem]{Corollary}
\theoremstyle{definition}
\newtheorem{defin}[theorem]{Definition}
\newtheorem{defins}[theorem]{Definitions}
\newtheorem{example}[theorem]{Example}
\newtheorem{xca}[theorem]{Exercise}

\newcommand{\bbb}{\mbox{$\beta$}}
\newcommand{\aaa}{\mbox{$\alpha$}}
\newcommand{\Aaa}{\mbox{$\mathcal A$}}
\newcommand{\ccc}{\mbox{$\mathcal C$}}
\newcommand{\ddd}{\mbox{$\delta$}}
\newcommand{\Ddd}{\mbox{$\Delta$}}
\newcommand{\eee}{\mbox{$\epsilon$}}
\newcommand{\Fff}{\mbox{$\mathcal F$}}
\newcommand{\Ggg}{\mbox{$\Gamma$}}
\renewcommand{\ggg}{\mbox{$\gamma$}}
\newcommand{\kkk}{\mbox{$\kappa$}}
\renewcommand{\lll}{\mbox{$\lambda$}}
\newcommand{\Lll}{\mbox{$\Lambda$}}
\newcommand{\rrr}{\mbox{$\rho$}}
\newcommand{\sss}{\mbox{$\sigma$}}
\newcommand{\Sss}{\mbox{$\mathcal S$}}
\newcommand{\Ss}{\mbox{$\Sigma$}}
\newcommand{\Th}{\mbox{$\Theta$}}
\newcommand{\ttt}{\mbox{$\tau$}}
\newcommand{\bdd}{\mbox{$\partial$}}
\newcommand{\zzz}{\mbox{$\zeta$}}
\newcommand{\inter}{\mbox{${\rm int}$}}

\title{On the additivity of knot width}

\authors{Martin Scharlemann\\Abigail Thompson}
\address{Mathematics Department, University of California\\Santa
Barbara, CA 93106, USA}
\secondaddress{Mathematics Department, University of California\\Davis, 
CA 95616, USA}
\asciiaddress{Mathematics Department, University of California\\Santa
Barbara, CA 93106, USA\\and\\Mathematics Department,
 University of California\\Davis, CA 95616, USA}
\gtemail{\mailto{mgscharl@math.ucsb.edu}, \mailto{thompson@math.ucdavis.edu}}
\asciiemail{mgscharl@math.ucsb.edu, thompson@math.ucdavis.edu}

\primaryclass{11Y16, 57M50}
\secondaryclass{57M25}
\keywords{Knot, width, additivity, Haken surfaces}

\begin{abstract}
It has been conjectured that the geometric invariant of knots in
$3$--space called the {\em width} is nearly additive.  That
is, letting $w(K) \in 2\mathbb{N}$ denote the width of a knot $K
\subset S^{3}$, the conjecture is that $w(K \# K') = w(K) +
w(K') - 2$.  We give an example of a knot $K_{1}$ so that for
$K_{2}$ any $2$--bridge knot, it appears that $w(K_{1} \# K_{2}) =
w(K_{1})$, contradicting the conjecture.
\end{abstract}

\asciiabstract{%
It has been conjectured that the geometric invariant of knots in
3-space called the width is nearly additive.  That is, letting w(K) in
N denote the width of a knot K in S^3, the conjecture is that w(K \#
K') = w(K) + w(K') - 2.  We give an example of a knot K_1 so that for
K_2 any 2-bridge knot, it appears that w(K_1 \# K_2) = w(K_1),
contradicting the conjecture.}

\maketitle
\cl{\small\it Dedicated to Andrew Casson, a mathematician's
mathematician}

\section{Background}

\begin{figure}[ht!]
\cl{
\relabelbox\small
\epsfxsize=0.5\textwidth
\epsfbox{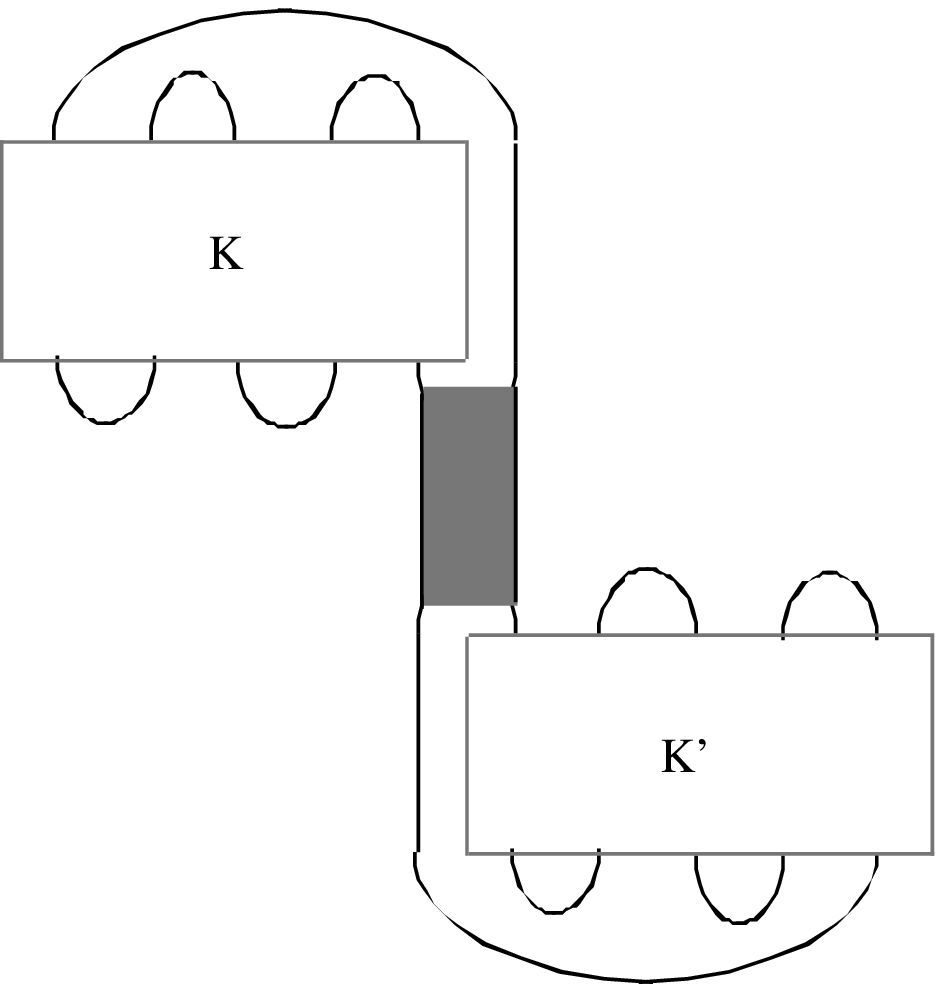}
\relabel{K}{$K$}
\relabel{K'}{$K'$}
\endrelabelbox}
\caption{}
\label{fig:widthsum}
\end{figure}

In \cite{G} Gabai associated to a knot in $3$--space an even number
called its width (see Definition \ref{def:width} below for the precise
definition).
Width can be viewed as a generalization of bridge number, an invariant
first studied by Schubert \cite{S}.  Schubert's remarkable discovery
was that, for $b(K)$ the bridge number of $K \subset S^{3}$, it is
always true that
$$b(K \# K') = b(K) + b(K') -1.$$
An alternative formulation of what Schubert showed is this: one way of
putting the knot sum $K \# K'$ into minimal bridge position is to put
both $K, K'$ into minimal bridge position, then place them vertically
adjacent in $S^{3}$.  Then create their knot sum via a vertical band
(see Figure \ref{fig:widthsum}).

Because width generalizes the notion of bridge number, it is natural
to hope that a similar equality is true for width.  Just as the
construction described above and shown in Figure \ref{fig:widthsum}
makes it obvious that
$$b(K \# K') \leq b(K) + b(K') - 1$$
(Schubert's deep contribution is the proof of the reverse
inequality, see \cite{Sch}), it also shows that knot width $w$
satisfies the inequality
$$w(K \# K') \leq w(K) + w(K') - 2.$$
That is, the resulting presentation of $K \# K'$ has width precisely
equal to $w(K) + w(K') - 2$.  The inequality reflects uncertainty over
whether this presentation is of minimal width (that is, {\em thin}) among
all possible presentations for $K \# K'$.  This suggests the following

\begin{conj} \label{conj:width}
For all knots $K, K' \subset S^{3}$, 
$$w(K \# K') = w(K) + w(K') - 2.$$
\end{conj}

Beyond the analogy with bridge number, and the ease of the inequality
in one direction, there are two pieces of evidence for Conjecture
\ref{conj:width}.  One piece of evidence is the proof in \cite{RS}
that the conjecture is true for knots that are {\em small}, that is,
knots which have no closed essential surfaces in their complements.  A
second but weaker piece of evidence is the proof in \cite{ScSc} that
in any case,
$$w(K \# K') \geq \max \{w(K), w(K') \}.$$
This shows that at least some inequality in the reverse direction is
true.

The aim here is to present an example of a knot (actually a family of
knots) $K_{1}$ which appears to have the property that
$$w(K_{1} \# K_{2}) = w(K_{1})$$
whenever $K_{2}$ is a $2$--bridge knot.  Since the width of a
$2$--bridge knot is $8$, this would be a counterexample to Conjecture
\ref{conj:width}.  Although we are not able to prove this equality,
the construction of $K_{1}$ is so flexible that the example does seem
to undercut any hope that the conjecture is true.  Only the absence so
far of a good method to prove that our presentation of the knot (family)
$K_{1}$ is the presentation of least width stands in the way of a complete
proof that at least some of these knots are counterexamples to Conjecture
\ref{conj:width}.

\section{The example}

\begin{figure}[ht!]
\cl{
\relabelbox\small
\epsfxsize=0.7\textwidth
\epsfbox{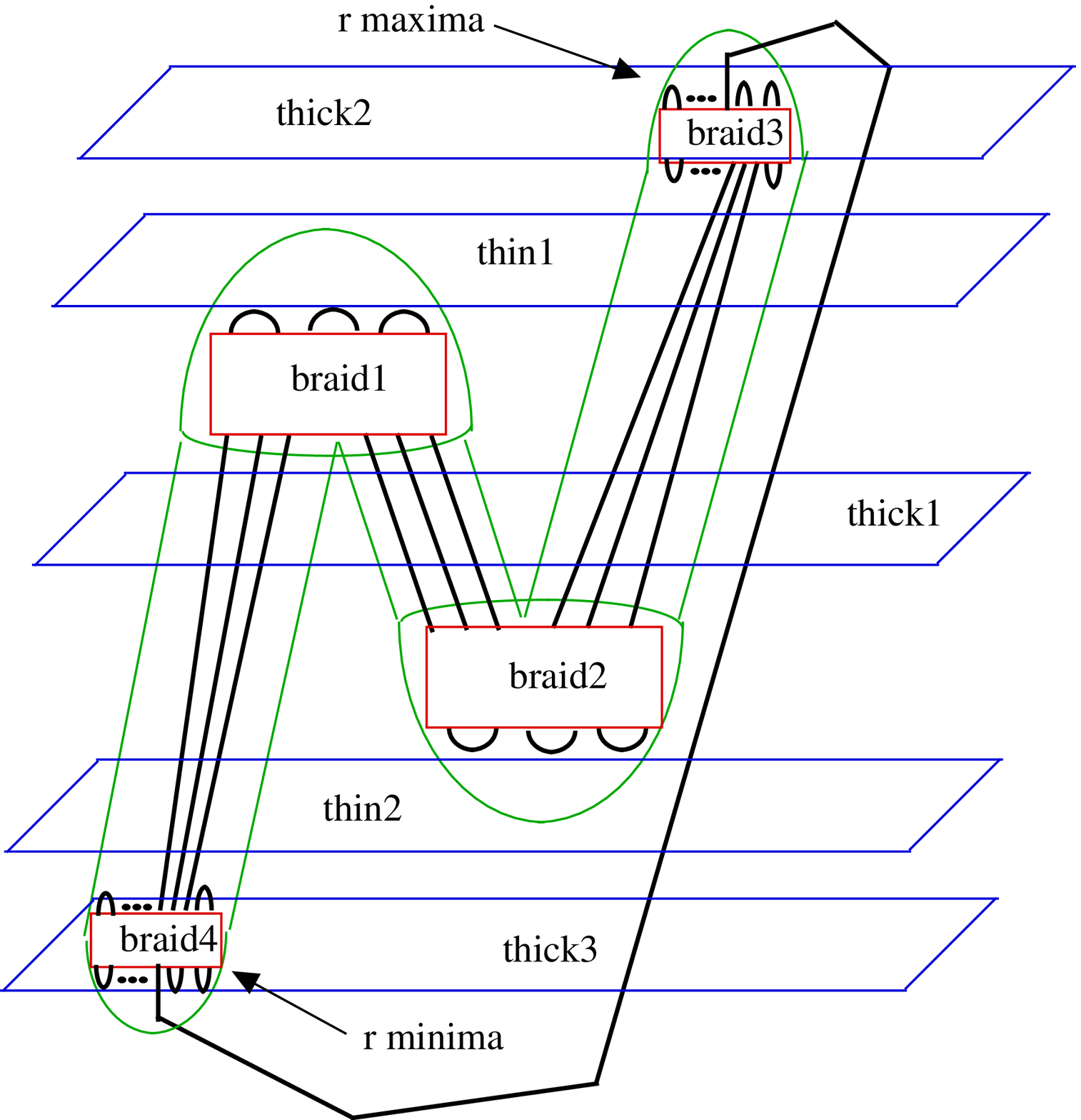}
\adjustrelabel <-3pt,-1pt> {braid1}{braid}
\adjustrelabel <-3pt,-1pt> {braid2}{braid}
\adjustrelabel <-3pt,-1pt> {braid3}{braid}
\adjustrelabel <-3pt,-1pt> {braid4}{braid}
\relabel{thin1}{thin}
\relabel{thin2}{thin}
\relabel{thick1}{thick}
\relabel{thick2}{thick}
\relabel{thick3}{thick}
\adjustrelabel <-6pt,3pt> {r maxima}{$r$ maxima}
\adjustrelabel <-3pt,0pt> {r minima}{$r$ minima}
\endrelabelbox
}
\caption{}
\label{fig:knotfromhell}
\end{figure}

The example $K_{1}$ is shown in Figure \ref{fig:knotfromhell}.  It
is actually a family of examples, because specific braids inside of
boxes are not specified; the point will be that, in search of an
example, the number of critical points $r$ can be made arbitrarily
high and the braids themselves be made arbitrarily complicated.

\begin{figure}[ht!]
\cl{
\epsfxsize=0.5\textwidth
\epsfbox{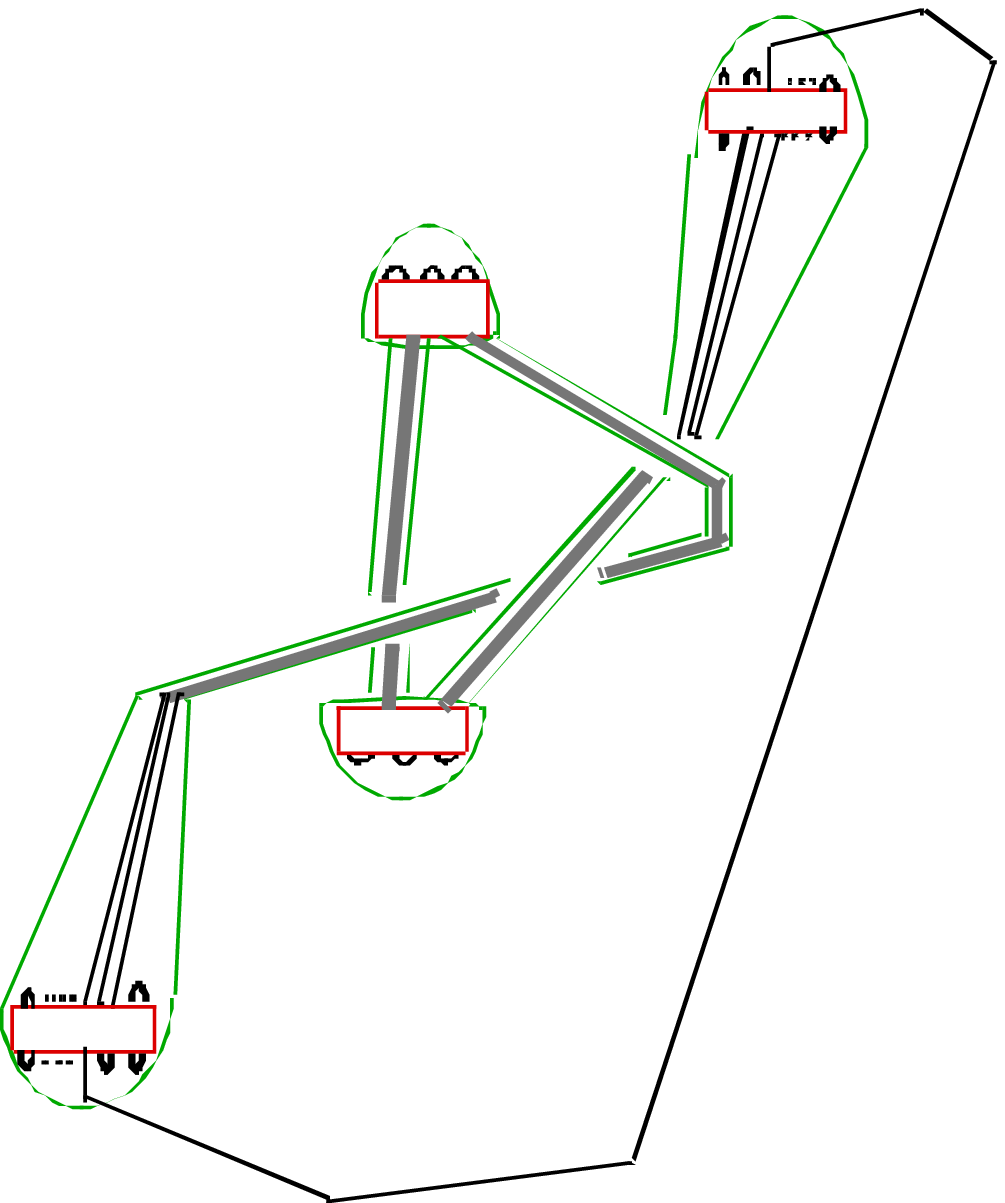}}
\caption{}
\label{fig:knotsumhell}
\end{figure}

Figure \ref{fig:knotsumhell} shows an imbedding of $K_{1} \# 
(\rm{trefoil})$ into
$S^{3}$; a decomposing sphere is shown in the figure.  Note that the
width of the presentations of $K_{1}$ and $K_{1} \# (\rm{trefoil})$
with respect to the vertical height function are the same, since
there is an obvious level-preserving reimbedding of one to the other.
(The construction of $K_{1}$ was inspired by the extensive use of
level-preserving reimbeddings in \cite{ScSc}.)  Moreover, the only
property of the trefoil knot that is used in the level-preserving
reimbedding is the fact that it is $2$--bridge --- any other $2$--bridge
knot would do.

\begin{figure}[ht!]
\cl{
\epsfxsize=0.7\textwidth
\epsfbox{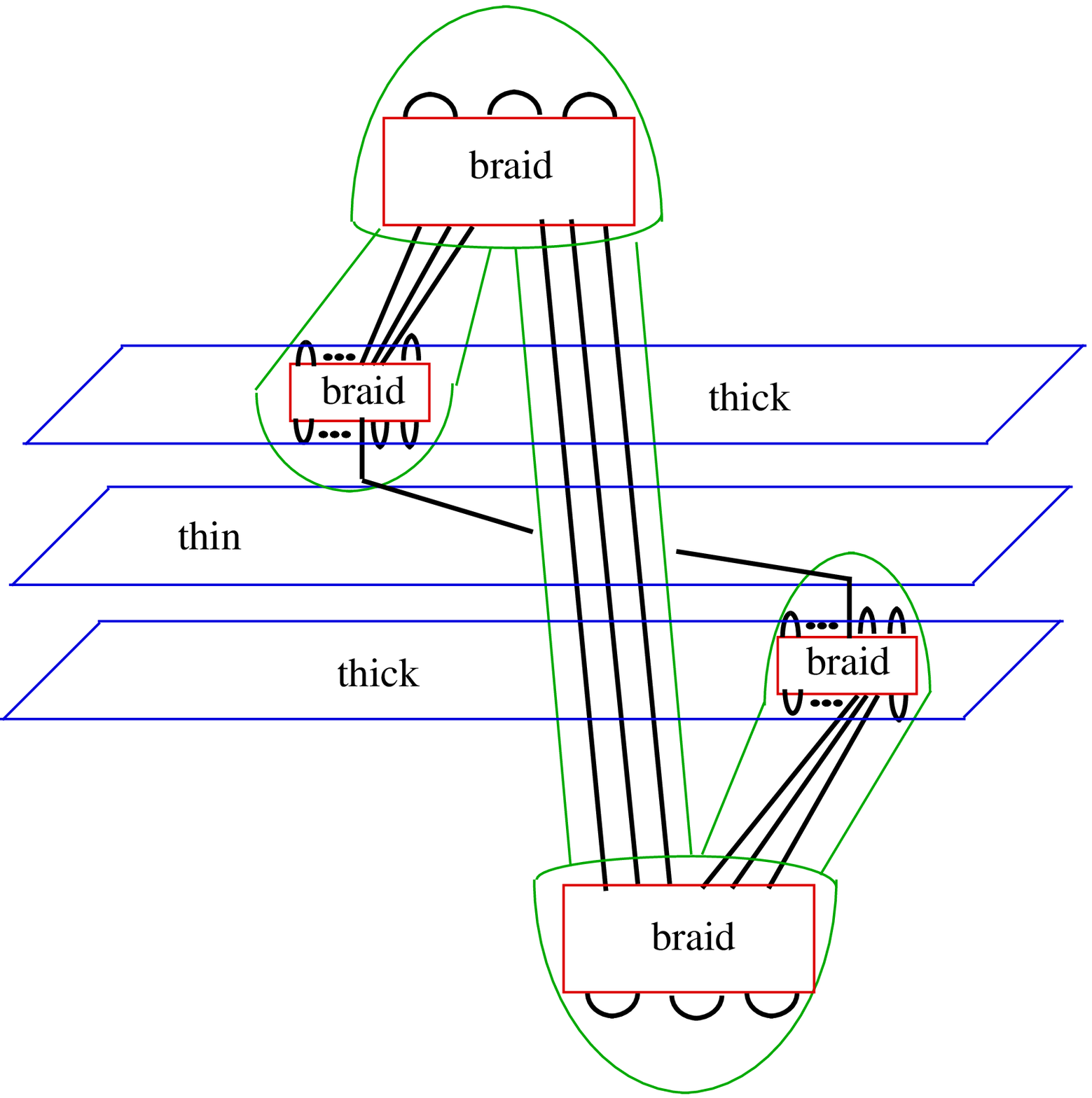}}
\caption{}
\label{fig:bridgefromhell}
\end{figure}

If we knew that the given presentation of the knot $K_{1}$ is thin (that
is, of minimal width) we would be done.  We do not yet know a full proof
of this, but cite two pieces of evidence:

First of all, the braids in the four braid boxes of Figure
\ref{fig:knotfromhell} may be made very complicated so that the number
of critical points of a plane sweeping across the boxes in any
direction except the vertical will have a large number of critical
points and arbitrarily high width.  Although this does not obviously
rule out the possibility of more complicated sweep-outs, it seems very
likely that with sufficiently complicated braids in place a sweep-out
can only be thin for the whole knot if the sweep-out passes through
these braid-boxes vertically, that is as horizontal planes, just as in
the given presentation.  (This is a stronger requirement than merely
ensuring that the tangles within the braid boxes are thin.)  Once this
is ensured, the order in which the sweep-out passes through each
tangle could be compared to the given order; alternative orders, such
as the order shown in Figure \ref{fig:bridgefromhell}, will give rise
to greater width (see calculation below) as long as $r \geq 4$.

This second but indirect bit of evidence that $K_{1}$ is thin is this:
According to Schubert's theorem, there is a presentation of $K_{1}$ with
lower bridge number than that of $K_{1} \# (\rm{trefoil})$ shown in Figure
\ref{fig:knotsumhell}.  But the Figure \ref{fig:knotsumhell}
presentation has the same bridge number as the presentation of $K_{1}$
shown in Figure \ref{fig:knotfromhell}.  Given that width is a
generalization of bridge number, we ought to check that the
presentation that lowers bridge number from that in Figure
\ref{fig:knotfromhell} does not also lower width.   Figure
\ref{fig:bridgefromhell} shows a presentation of $K_{1}$ that has 
lower bridge number than the original presentation, shown in Figure 
\ref{fig:knotfromhell}.  But in fact the presentation
of $K_{1}$ in Figure \ref{fig:bridgefromhell} is not as thin as that in Figure
\ref{fig:knotfromhell}.

To explicitly compare the widths of the two presentations, recall the
definitions (see \cite{ScSc}):

\begin{defin}
\label{def:width}
Let $p\co S^3 \rightarrow \mathbb{R}$ be the standard height function
and for each
$-1 < t < 1$ let $S^{t}$ denote the sphere $p^{-1}(t)$.  Let $K
\subset S^3$ be a knot in general position with respect to $p$ and
$c_1, \dots, c_n$ be the critical values of $p|K$ listed in
increasing order $c_1 < \dots < c_n$.  Choose $r_1, \dots, r_{n-1}$ so
that $c_i < r_i < c_{i+1}$ for $i = 1, \ldots, n-1$.  The {\em width of K
with respect to $p$}, denoted by $w(K, p)$, is $\sum_i|K \cap
S^{r_{i}}|$.  The {\em width of K}, denoted by $w(K)$, is the minimum
of $w(K', p)$ over all knots $K'$ isotopic to $K$.  We say that $K$ is
in {\em thin position} if $w(K, p) = w(K).$
\end{defin}

The key operational feature of width is that it is reduced if a
minimum is pushed up above a maximum whereas switching the levels of
two adjacent minima or two adjacent maxima has no effect.

There is an easier way to calculate width.  For the levels $r_{i}$
described above, call $r_i$ a {\em thin level} of $K$ (and the sphere 
$S^{r_{i}}$ a {\em thin sphere} for $K$) with respect to
$p$ if $c_i$ is a maximum value for $p|K$ and $c_{i+1}$ is a minimum
value for $p|K$.  Dually $r_i$ is a {\em thick level} of $K$ 
(and the sphere 
$S^{r_{i}}$ a {\em thick sphere} for $K$) with
respect to $p$ if $c_i$ is a minimum value for $p|K$ and $c_{i+1}$ is a
maximum value for $p|K$.  Since the lowest critical point of $p|K$ is a
minimum and the highest is a maximum, there is one more thick level
than thin level.

\begin{lemma}
Let $r_{i_1}, \dots, r_{i_k}$ be the thick levels of $K$ and $r_{j_1},
\dots, r_{j_{k-1}}$ the thin levels.  Set $a_{l} = |\;K\; \cap
S^{r_{i_l}}|$ and $b_{l} = |\;K\; \cap S^{r_{j_l}}|$.  Then
$$2w(K) = \sum_{l = 1}^k a_{l}^2 -  \sum_{l = 1}^{k-1} b_{l}^2.$$
\end{lemma}
\begin{proof} See \cite{ScSc}. \end{proof}

Now apply this formula to the presentations of $K_{1}$ given in Figures 
\ref{fig:knotfromhell} and \ref{fig:bridgefromhell}, with thick and 
thin spheres (appearing as planes) noted in the figures.  The former has width 
$$4(r+1)^{2} + 50 - 16 = 2(2r^{2} + 4r + 19).$$
The latter has width
$$4(r+2)^{2} - 8 = 2(2r^{2} + 8r + 4).$$
So as long as $r \geq 4 \Rightarrow 8r + 4 > 4r + 19$ the presentation
of $K_{1}$ in Figure \ref{fig:knotfromhell} is thinner.

\section{Higher bridge number}

It is reasonable to ask whether likely counterexamples to the
conjecture are limited to knots, such as those above, whose width is
unchanged by adding $2$--bridge knots.  After all, $2$--bridge knots
often play a special role in knot theory.  Almost certainly the answer
is no; in this section we briefly note how a $3$--bridge counterexample
might be constructed. In principle the same ideas should
work with arbitrary bridge number, though this would appear to be
increasingly difficult to demonstrate.

\begin{figure}[ht!]
\cl{
\relabelbox\small
\epsfxsize=0.5\textwidth
\epsfbox{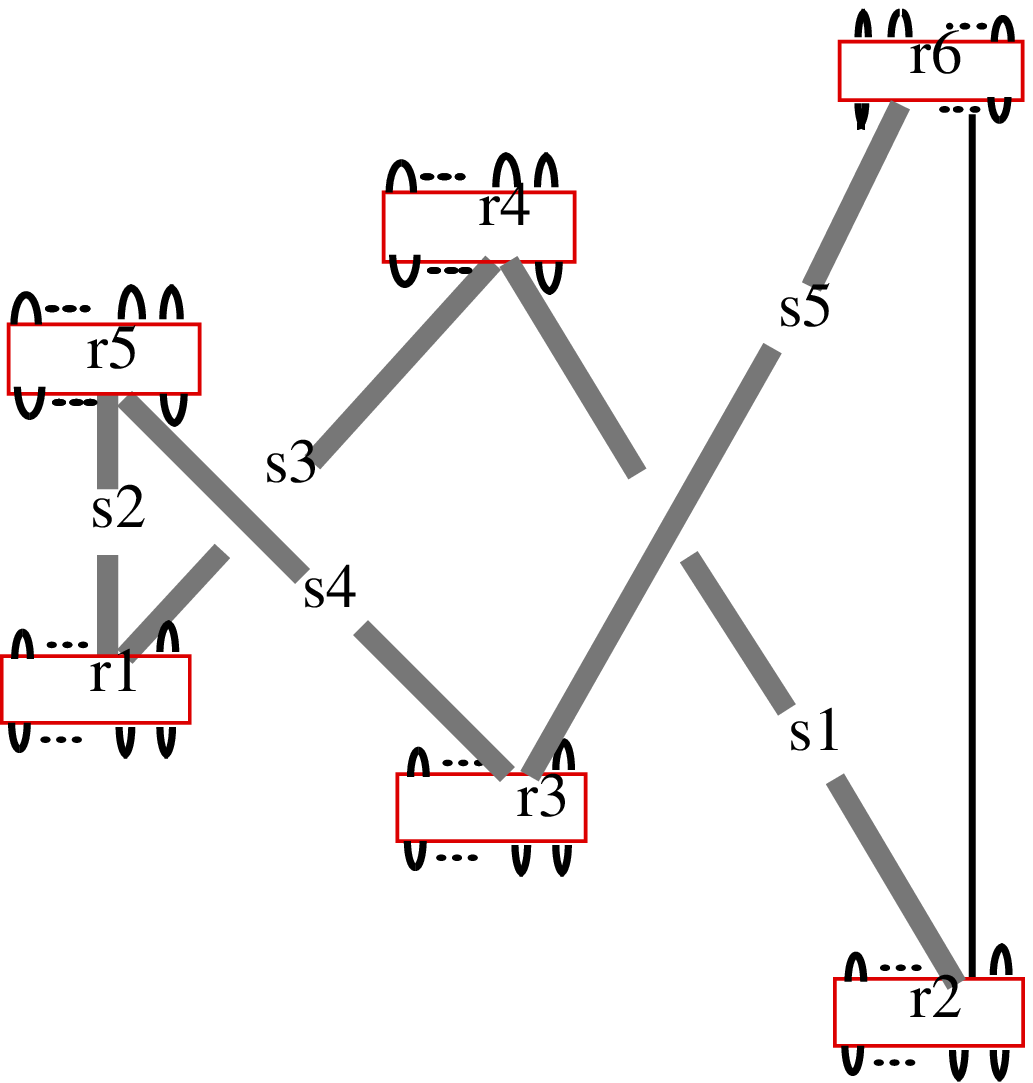}
\adjustrelabel <-4pt,-2pt> {r1}{$r_1$}
\adjustrelabel <-2pt,-2pt> {r2}{$r_3$}
\adjustrelabel <-8pt,-1pt> {r3}{$r_2$}
\adjustrelabel <-4pt,-2pt> {r4}{$r_2$}
\adjustrelabel <-2pt, 0pt> {r5}{$r_1$}
\adjustrelabel < 0pt,-2pt> {r6}{$r_3$}
\adjustrelabel <-2pt, 0pt> {s1}{$s_3$}
\adjustrelabel <-2pt,-1pt> {s2}{$s_1$}
\adjustrelabel <-3pt,-3pt> {s3}{$s_2$}
\adjustrelabel < 0pt, 0pt> {s4}{$s_2$}
\adjustrelabel <-2pt, 0pt> {s5}{$s_3$}
\endrelabelbox
}
\caption{} \label{fig:3bridge}
\end{figure}  

Consider a knot $K_{3}$ of the form described in Figure
\ref{fig:3bridge}.  The notation $r_{i}$ denotes the number of strands
of the part of the knot that lies inside the ``braid-box'', that is,
twice the number of maxima appearing above the upper boxes (respectively,
below the lower boxes).  The notation $s_{i}$ refers to the number of
parallel strands between the boxes.  The thin unmarked strand
represents a single strand.  We have earlier seen that the knot
$K_{1}$ in Figure \ref{fig:knotfromhell} can be added to any
$2$--bridge knot, without apparently affecting its width, by
reimbedding it in a level-preserving way as a satellite of any
$2$--bridge knot.  In the same way, the knot $K_{3}$ in Figure
\ref{fig:3bridge} can be reimbedded in a level-preserving way as a
satellite of any $3$--bridge knot $K_{4}$.  This reimbedding gives a
presentation of the sum $K_{3} \# K_{4}$ that has the same width as the
original presentation of $K_{3}$.  So
if the given presentation of $K_{3}$ is thin, we have a counterexample
to Conjecture \ref{conj:width} in which one of the summand knots is 
the specific knot $K_{3}$ and the other is any $3$--bridge knot.

\begin{figure}[ht!]
\cl{
\relabelbox\small
\epsfxsize=0.5\textwidth
\epsfbox{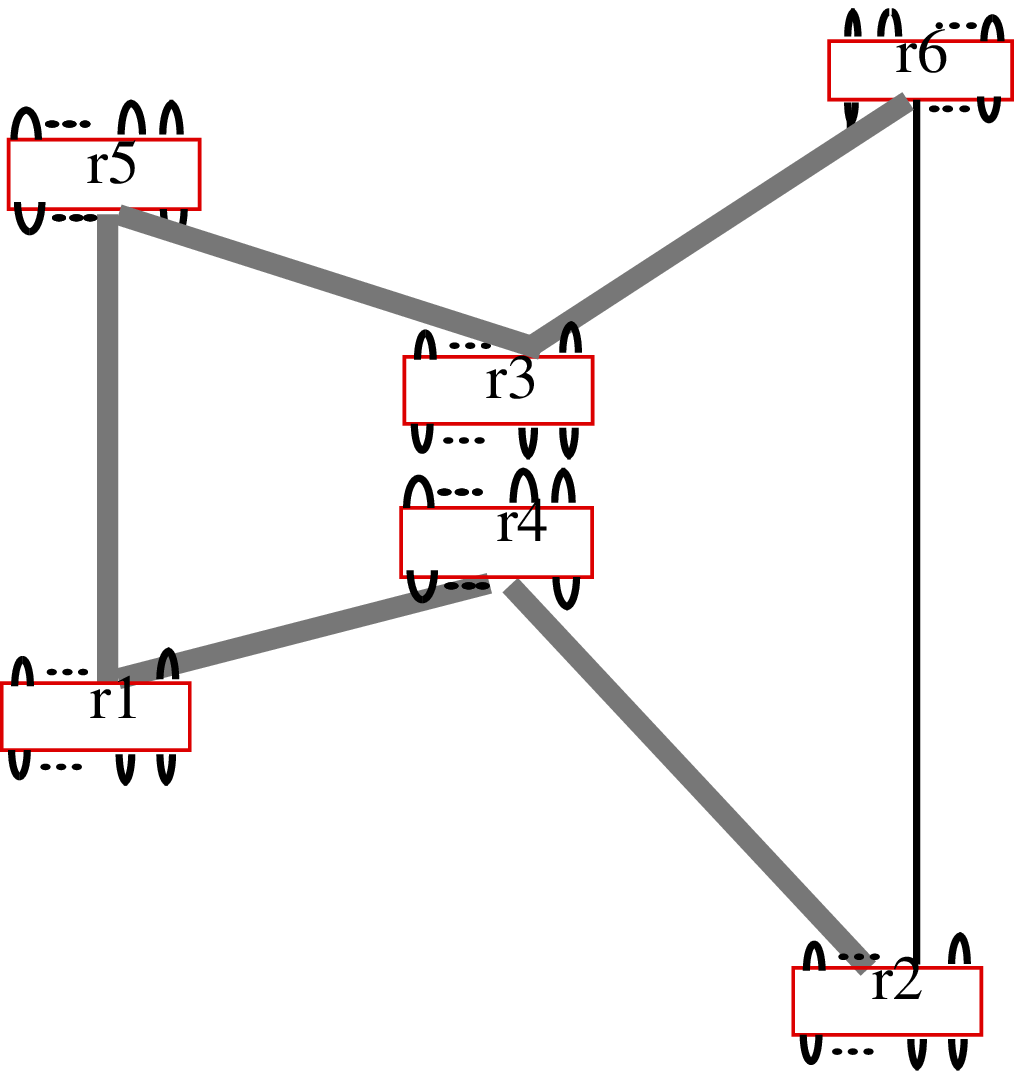}
\adjustrelabel <-2pt,-2pt> {r1}{$r_1$}
\adjustrelabel < 0pt,-2pt> {r2}{$r_3$}
\adjustrelabel <-2pt, 0pt> {r3}{$r_2$}
\adjustrelabel <-4pt,-2pt> {r4}{$r_2$}
\adjustrelabel <-2pt, 0pt> {r5}{$r_1$}
\adjustrelabel < 0pt,-2pt> {r6}{$r_3$}
\endrelabelbox
}
\caption{} \label{fig:3bridge2}
\end{figure}

Of course there is still the difficulty of showing that the
presentation of $K_{3}$ in Figure \ref{fig:3bridge} is thin.  But again the
flexibility in how the braid-boxes are filled in suggests that filling
in with complicated braids will force any thin presentation also to
sweep across the braid boxes vertically.  So a thin presentation is
likely to preserve the braid-box structure.  It is perhaps less
plausible here than in the case of $K_{1}$ that a presentation which
preserves the braid-box structure cannot be thinner than the original
presentation; indeed, an example of a significantly different
presentation of $K_{3}$ that preserves the braid-box structure is
given in Figure \ref{fig:3bridge2}.

However, it is straightforward (though not particularly easy) to show that with
appropriate choice of $r_{i}, s_{j}$ the first presentation is thinner
than the second.  A comparison of the two widths is unaffected by the
size of $r_{3}$.  In addition, so long as we

\begin{itemize}
\item take the other values of $r_{i}, s_{j}$ sufficiently large
\item set $r_{1}$ only minimally larger than $s_{1} + s_{2}$
\item set $s_{1} > s_{3}$, and, with these numbers set, 
\item take $r_{2} > s_{2} + s_{3}$ sufficiently high,
\end{itemize}
then the second presentation of $K_{3}$ will be
wider than the first.

\begin{figure}[ht!]
\cl{
\relabelbox\small
\epsfxsize=0.5\textwidth
\epsfbox{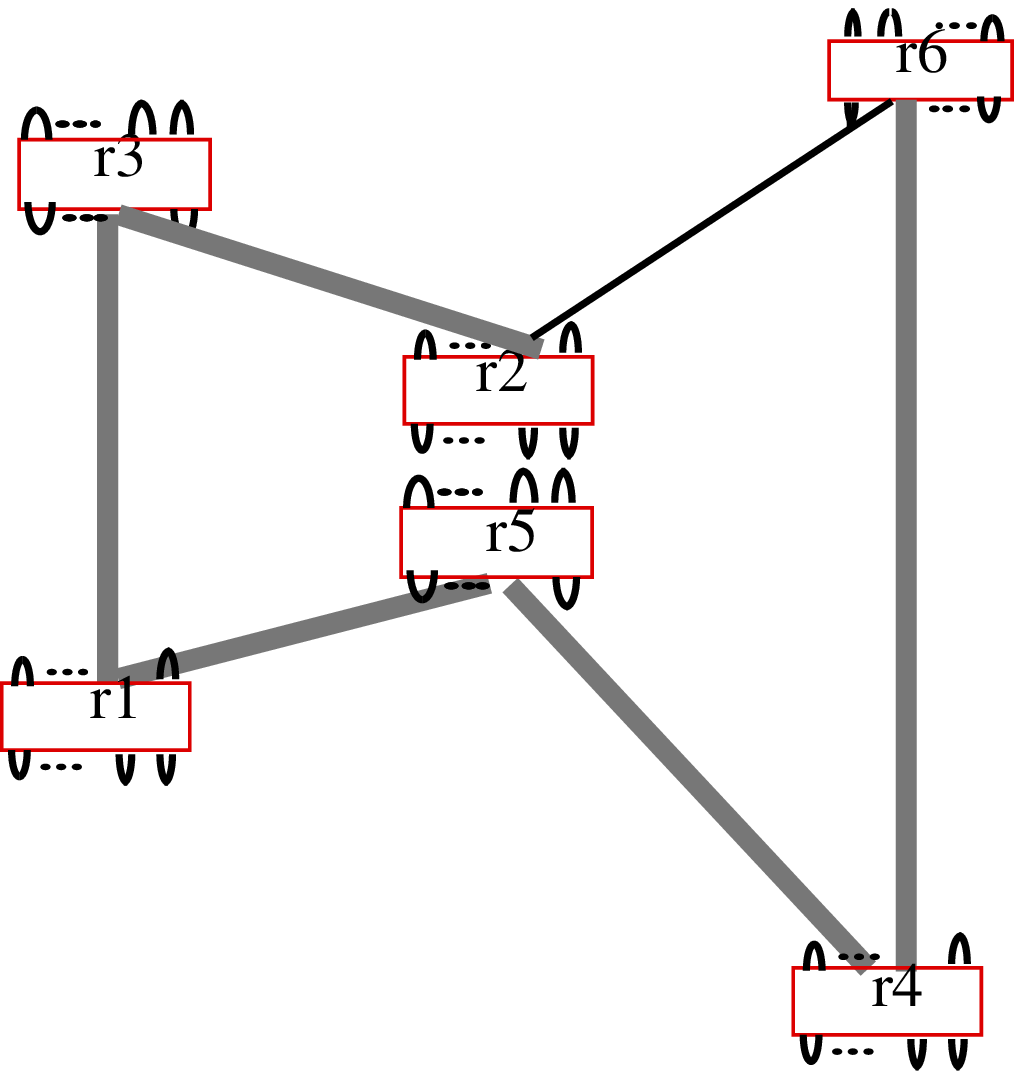}
\adjustrelabel <-2pt,-2pt> {r1}{$r_1$}
\adjustrelabel < 0pt,-2pt> {r2}{$r_3$}
\adjustrelabel <-2pt,-1pt> {r3}{$r_2$}
\adjustrelabel < 0pt,-1pt> {r4}{$r_2$}
\adjustrelabel <-2pt, 0pt> {r5}{$r_1$}
\adjustrelabel < 0pt,-2pt> {r6}{$r_3$}
\endrelabelbox
}
\caption{} \label{fig:3bridge3}
\end{figure}

Another typical braid-box preserving presentation of $K_{3}$ is given
in Figure \ref{fig:3bridge3}.  For this presentation, width can be
made arbitrarily high simply by increasing $r_{3}$ sufficiently.  As
already noted, this has no effect on the comparison above, so
ultimately we have values for $r_{i}, s_{j}$ that guarantee the first
presentation is thinnest of the three.  Other braid-box preserving
presentations are easily dealt with in a similar manner.  As was the
case for the $2$--bridge examples of the previous section, there
remains the challenge of showing that putting sufficiently complicated
braids into the braid-boxes will guarantee that a thin presentation
will preserve those braid-boxes.

\section*{Acknowledgement}
This research was partially supported by NSF grants.

\Addresses\recd

\end{document}

%% file: gtmonout.tex
\def\ifplaintex{\expandafter\ifx\csname documentclass\endcsname\relax}

\ifplaintex 
\else
\fi

\def\gt{{\mathsurround=0pt\it $\cal G\mskip-2mu$eometry \&\ 
$\cal T\!\!$opology}}        

\def\gtm{{\mathsurround=0pt\it $\cal G\mskip-2mu$eometry \&\ 
$\cal T\!\!$opology $\cal M\mskip-1mu$onographs}}    

\def\gtp{{\mathsurround=0pt\it $\cal G\mskip-2mu$eometry \&\ 
$\cal T\!\!$opology $\cal P\!$ublications}}  

\def\recd{{\small Received:\qua\receiveddate\ifx\reviseddate\relax
\else\qquad Revised:\qua\reviseddate\fi\par}} 

\def\volumenumber#1{\def\thevolumenumber{#1}}
\def\volumeyear#1{\def\thevolumeyear{#1}}
\def\volumename#1{\def\thevolumename{#1}}
\def\papernumber#1{\def\thepapernumber{#1}}
\def\pagenumbers#1#2{\def\startpage{#1}\def\finishpage{#2}}
\def\published#1{\def\publishdate{#1}}
\def\received#1{\def\receiveddate{#1}}
\def\revised#1{\def\reviseddate{#1}}
\def\accepted#1{\def\accepteddate{#1}}

\def\asciiaddress#1{\def\theasciiaddress{#1}}
\def\asciiemail#1{\def\theasciiemail{#1}}

\long\def\asciiabstract#1{\long\def\theasciiabstract{#1}}

\let\\\par
\let\thevolumenumber\relax\let\thepapernumber\relax
\let\thevolumeyear\relax\let\startpage\relax
\let\finishpage\relax\let\publishdate\relax\let\receiveddate\relax
\let\reviseddate\relax\let\accepteddate\relax\let\theasciititle\relax
\let\theasciiauthors\relax\let\theasciiaddress\relax
\let\theasciiabstract\relax

\let\theerratum\relax\let\theasciiemail\relax
\let\theshortauthors\relax\let\theshorttitle\relax

\def\startpage{1}\def\finishpage{15}\def\thepapernumber{77}

\volumenumber{2}
\volumename{Proceedings of the Kirbyfest}
\volumeyear{1999}

\long\def\maketitlep{   

\count0=\startpage

\gtm\nl        
{\small Volume \thevolumenumber: \thevolumename\nl 
\ifx\theerratum\relax\else Erratum \erratumnumber\nl\fi
Pages \startpage--\finishpage\nl}

\vglue 0.1truein   

{\parskip=0pt\leftskip 0pt plus 1fil\def\\{\par\smallskip}{\ifplaintex\large
\else\Large\fi\bf\thetitle}\par\medskip}   
\vglue 0.05truein 

{\parskip=0pt\leftskip 0pt plus 1fil\def\\{\par}{\sc\theauthors}
\par\medskip}%
 
\vglue 0.03truein 

{\small\leftskip 25pt\rightskip 25pt{\bf Abstract}\stdspace\theabstract

{\bf AMS Classification}\stdspace\theprimaryclass
\ifx\thesecondaryclass\relax\else; \thesecondaryclass\fi\par
{\bf Keywords}\stdspace \thekeywords\par}\vglue 7pt

}   

\font\phead=cmsl9 scaled 950
\font=cmsl9 scaled 1050
\font\pnum=cmbx10 scaled 913
\font\lnum=cmbx10 
\font\pfoot=cmsl9 scaled 950
\font=cmsl9 scaled 1050
\ifplaintex
\headline{\vbox to 0pt{\vskip -4.5mm\line{\small\phead\ifnum
\count0=\startpage ISSN 1464-8997 (on line)
1464-8989 (printed) \hfill {\pnum\folio}\else\ifodd\count0\def\\{ }%
\ifx\theshorttitle\relax\thetitle\else\theshorttitle\fi\hfill{\pnum\folio}
\else\def\\{ and }{\pnum\folio}\hfill\ifx\theshortauthors\relax\theauthors
\else\theshortauthors\fi\fi\fi}\vss}}
\footline{\vbox to 0pt{\vglue 0mm\line{\small\pfoot\ifnum\count0=\startpage
Published \publishdate:\qua\copyright\ \gtp\hfill\else
\gtm, Volume \thevolumenumber\ (\thevolumeyear)\hfill\fi}\vss
}}
\else
\makeatletter
\def\@oddhead{{\small\lhead\ifnum\count0=\startpage ISSN 1464-8997 (on line)
1464-8989 (printed) \hfill {\lnum\number\count0}\else\ifodd\count0
\def\\{ }\ifx\theshorttitle\relax \thetitle \else\theshorttitle\fi\hfill
{\lnum\number\count0}\else\def\\{ and }{\lnum\number\count0}
\hfill\ifx\theshortauthors\relax 
\theauthors\else\theshortauthors\fi\fi\fi}}\def\@evenhead{@oddhead}
\def\@oddfoot{\small\lfoot\ifnum\count0=\startpage Published \publishdate:\qua\copyright\ \gtp\hfill\else
\gtm, Volume \thevolumenumber\ (\thevolumeyear)\hfill\fi}
\def\@evenfoot{@oddfoot}
\makeatother
\fi

\let\maketitlepage\maketitlep
\let\makeshorttitle\maketitlepage
\let\maketitle\maketitlepage

\newwrite\gtoutfile
\long\gdef\makeheadfile{  
{\def\\{, }\def\s{ }
\immediate\openout\gtoutfile head.xxx
\immediate\write\gtoutfile{Proxy-for: \ifx\theasciiauthors\relax
\theauthors\else\theasciiauthors\fi\s<\ifx\theasciiemail\relax\theemail\else\theasciiemail\fi>}
\immediate\write\gtoutfile{\noexpand\\}
\immediate\write\gtoutfile{Authors: \ifx\theasciiauthors\relax
\theauthors\else\theasciiauthors\fi}
{\def\\{ }\immediate\write\gtoutfile{Title: \ifx\theasciititle\relax
\thetitle\else\theasciititle\fi}}
\immediate\write\gtoutfile{Subj-class: GT or SG, GR etc}
\immediate\write\gtoutfile{MSC-class: \theprimaryclass\ifx\thesecondaryclass\relax\else, \thesecondaryclass\fi}
\immediate\write\gtoutfile{Journal-ref: Geom. Topol. Monogr. \thevolumenumber\s
(\thevolumeyear) \startpage-\finishpage}
\immediate\write\gtoutfile{Comments: Published by Geometry and Topology Monographs at}
\immediate\write\gtoutfile{\s\s\s  http://www.maths.warwick.ac.uk/gt/GTMon\thevolumenumber/paper\thepapernumber.abs.html}
\immediate\write\gtoutfile{\noexpand\\}
\immediate\write\gtoutfile{}
\ifx\theasciiabstract\relax
\immediate\write\gtoutfile{\theabstract}\else
\immediate\write\gtoutfile{\theasciiabstract}\fi
\immediate\write\gtoutfile{}
\immediate\write\gtoutfile{\noexpand\\}
\immediate\write\gtoutfile{}
\immediate\closeout\gtoutfile}}  

\def\maketitlepage{\maketitlep\makeheadfile}
\let\makeshorttitle\maketitlepage
\let\maketitle\maketitlepage